\documentclass[a4paper]{amsproc}

\usepackage[resetstyle]{phfthm}
\usepackage{amssymb}
\usepackage{pstricks}
\usepackage{pstcol}
\usepackage{graphicx}
\usepackage{amsmath}
\usepackage{amsfonts}
\usepackage{latexsym}
\usepackage[urlcolor=blue, colorlinks=true, linkcolor=blue, citecolor=blue]{hyperref}
\theoremstyle{plain}

\theoremstyle{definition}

\numberwithin{equation}{section}
\theoremstyle{definition}
\newtheorem{ex}{Example}
\renewcommand{\leq}{\leqslant}
\renewcommand{\geq}{\geqslant}
\renewcommand{\setminus}{\smallsetminus}

\DeclareMathOperator{\li}{Li}

\newcommand{\rmd}{\mathrm}

\newcommand{\G}{\textbf{\textup{G}}}

\newcommand{\C}{{\mathbb C}}
\newcommand{\N}{{\mathbb N}}

\title[New closed forms for a dilogarithmic integral, related integrals, and series]{New closed forms for a dilogarithmic integral, related integrals, and series}



\author[Abdulhafeez A. Abdulsalam]{Abdulhafeez A. Abdulsalam}

\address{\hfill{\it Received 03 05
2024, revised 10 06 2024}\newline Department of Mathematics, University of Ibadan, \newline Ibadan, Oyo, \newline
Nigeria
}

\email{hafeez147258369@gmail.com}

\allowdisplaybreaks
\begin{document}

{\begin{flushleft}\baselineskip9pt\scriptsize {\bf SCIENTIA}\newline
Series A: {\it Mathematical Sciences}, Vol. 35 (2025), 1--xx
\newline Universidad T\'ecnica Federico Santa Mar{\'\i}a
\newline Valpara{\'\i}so, Chile
\newline ISSN 0716-8446
\newline {\copyright\space Universidad T\'ecnica Federico Santa
Mar{\'\i}a\space 2025}
\end{flushleft}}

\vspace{10mm} \setcounter{page}{1} \thispagestyle{empty}

\begin{abstract}
In this study, we present a new closed form for the generalized integral 
$$\int_0^1 \frac{\li_2(z) \ln(1+az)}{z}\, \rmd{d}z,$$
where $a \in \C \setminus(-\infty, -1)$ and $\li_2(z)$ is the dilogarithm function. This generalization is achieved by leveraging our established findings in conjunction with V\u{a}lean's results. Furthermore, we provide explicit closed forms for associated integrals, prove a transformation formula for double infinite series, expressing them as the sum of the square of an infinite series and another infinite series. We utilize this relationship to derive a novel closed form for the generalized series
$$\sum_{k=1}^\infty \frac{ \zeta\left(m, \frac{rk-s}{r}\right) }{(rk-s)^m},$$
for $\Re(m) > 1$, $r, s \in \mathbb{C}$, where $r \neq 0$, $rk \neq s$, for any positive integer $k$, and $\zeta(s, z)$ denotes the Hurwitz zeta function. Utilizing Hermite's integral representation for $\zeta(s, z)$, we derive a family of integrals from this series.
\end{abstract}

\maketitle

\section{Introduction}
In this paper, we provide a new closed form for the integral 
\begin{equation}\label{ruh1}
\int_0^1 \frac{\li_2(z) \ln(1+az)}{z}\, \rmd{d}z,
\end{equation}
where $\li_2(z)$ denotes the dilogarithm function, defined as \cite[(25.12.1)]{bib23}
\[\li_2(z) = \sum_{k=1}^\infty \frac{z^k}{k^2}, \quad |z| \leq 1.\]
Our approach entails the transformation of integrals into infinite series involving harmonic numbers, followed by the subsequent evaluation of these resulting series. Through this method, we not only determine a closed form for the aforementioned integral but also discover a closed form for
\[\sum_{k=1}^\infty \frac{ \zeta\left(m, \frac{rk+r-s}{r}\right) }{(rk-s)^m}, \quad \Re{m} > 1\wedge\hspace{0.03cm} r, s \in \C, r \neq 0, rk \neq s, \forall\hspace{0.1cm}  k\, \in \N.\]
 Utilizing Hermite's integral representation for $\zeta(s, z)$, we derive a family of integrals from this series. Examples of such integrals are:
\begin{align*}
&\int_0^\infty \frac{x}{(x^2 + 1)^2 (e^{2\pi x} - 1)} \, \rmd{d}x + \frac{1}{2} \int_0^\infty \frac{x}{(x^2 + 4)^2 (e^{2\pi x} - 1)} \, \rmd{d}x 
\\&\qquad+  \frac{1}{3} \int_0^\infty \frac{x}{(x^2 + 9)^2 (e^{2\pi x} - 1)} \, \rmd{d}x  +  \frac{1}{4} \int_0^\infty \frac{x}{(x^2 + 16)^2 (e^{2\pi x} - 1)} \, \rmd{d}x  + \cdots 
\\&=  \frac{\pi^4}{288} - \frac{\zeta(3)}{4},
\end{align*}
\begin{align*}
&\int_0^\infty \frac{x}{(4 x^2 + 1)^2 (e^{2\pi x} - 1)} \, \rmd{d}x + \frac{1}{3} \int_0^\infty \frac{x}{(4 x^2 + 9)^2 (e^{2\pi x} - 1)} \, \rmd{d}x 
\\&\qquad+  \frac{1}{5} \int_0^\infty \frac{x}{(4 x^2 + 25)^2 (e^{2\pi x} - 1)} \, \rmd{d}x  +  \frac{1}{7} \int_0^\infty \frac{x}{(4 x^2 + 49)^2 (e^{2\pi x} - 1)} \, \rmd{d}x  + \cdots 
\\&= \frac{\pi ^4}{1024}-\frac{7 \zeta (3)}{128}.
\end{align*}
These integrals do not appear in existing literature. Throughout this work, $H_n$ represents the $n$th harmonic number, defined as 
\[H_n = \sum_{k=1}^n \frac{1}{k}, \quad n \in \N,\]
$\psi_{m-1}(z)$ represents the polygamma function, defined as \cite[\S5.15]{bib23},
\[\psi_{m-1}(z) = (-1)^{m} (m - 1)!\sum_{k=0}^\infty \frac{1}{(k+z)^m}, \quad m \geq 2,\, m \in \N, z \not\in -\mathbb{N}_0,\] 
and $\zeta(s, z)$ represents the Hurwitz zeta function, defined as \cite[\S 25.11]{bib23}
\[\zeta(s, z) = \sum_{n=0}^{\infty} \frac{1}{(n+z)^s}, \quad z \not\in -\mathbb{N}_0, \Re\,s > 1.\]
By  incorporating V\u{a}lean's closed form for 
\[\int_0^1 \frac{\li_2(z) \ln(z)}{1-az}\, \rmd{d}z, \quad a \in \C \setminus(1, \infty) \cup \{0\},\]
alongside our derived closed forms for
\begin{equation}\label{begeqhh}
\int_0^1 \frac{\ln{z}\ln(1+az)\ln(1-z)}{z} \, \rmd{d}z, \quad \int_0^1 \frac{\ln{z}\ln^2(1 + a z)}{z}\, \rmd{d}z, \quad a \in \C \setminus(-\infty, -1),
\end{equation}
we present a generalized version of \eqref{ruh1} in the form  
\begin{equation}\label{begeqh1}
\int_0^1 \frac{\li_2(z) \ln(1+az)}{z}\, \rmd{d}z, \quad a \in \C \setminus(-\infty, -1).
\end{equation}
The two integrals in \eqref{begeqhh} are equal when $a=-1$. The simplest evaluation of \eqref{begeqh1} occurs when $a=-1$. In this case, we have
\begin{equation}
\int_0^1 \frac{\li_2(z) \ln(1-z)}{z}\, \rmd{d}z = -\frac{\pi^4}{72}.
\end{equation}
The exclusion of \eqref{begeqh1} for $a \in (-\infty, -1)$ stems from the observation that if $a \in (-\infty, -1)$, the integral diverges at a certain point within the integration domain. To illustrate this, set $a=-b$, where $b \in (1, \infty)$, and notice that for all values of $b$ in this interval, there exists a unique $z = \frac{1}{b} \in (0, 1)$ such that
\begin{equation}\label{beghy1}
 \frac{\li_2(z) \ln(1-bz)}{z}  \to -\infty.
\end{equation}
In summary, the results established in this article are outlined as follows. In Section \ref{sec3.1}, we present a closed form for \eqref{ruh1} using known results. We also establish the relationship $2\sum_{k=1}^\infty \frac{H_{2k}}{(2k)^3} = \int_0^1 \frac{\li_2(z) \ln(1+z)}{z}\, \rmd{d}z$. We begin our exploration by introducing a novel generalization of \eqref{ruh1}. Additionally, we provide generalizations for related integrals. Furthermore, we offer a representation for the series $\sum_{k=1}^\infty \frac{(-1)^{k} H_k a^k}{k^3}$ that allows us to provide the generalization \eqref{begeqh1} for $a \in \C\setminus(-\infty, 0]$ while avoiding logarithm of negative real numbers. Theorems \ref{anathm} and \ref{newthgh} are not new, as we use our established results to provide a new proof of Jonqui\`ere's inversion formula for order 4 and arguments $-\frac{1}{z}$ and $\frac{z}{z-1}$. The closed forms  presented in Theorems \ref{thmabd1}--\ref{thmhaf5} and \ref{theoremh9}--\ref{bghaft} are new and have not been presented elsewhere in the literature. In Section  \ref{sec3.2}, we introduce a transformative approach for double infinite series, enabling us to express them as sums of the square of an infinite series and another infinite series. We apply this theorem to derive novel generalized identities. The Computer Algebra System (CAS) software employed for result verification throughout this paper is \textsf{Mathematica 13}.
\section{Notations and Definitions}\label{sec2}
In this manuscript, we employ the following abbreviated notations: $B_n$ represents the Bernoulli numbers \cite[\S24.2(i)]{bib23}, $E_n$ represents the Euler numbers \cite[\S24.2(ii)]{bib23}, $\gamma \approx 0.5772156649$ represents Euler's constant, $\G = \sum_{k=0}^\infty \frac{(-1)^k}{(2k+1)^2}$ represents Catalan's constant ($\G \approx 0.9159655941$), while $\mathrm{e} \approx 2.71828182845$ stands for Euler's number. We define $\mathbb{N}_0 := \mathbb{N} \cup \{0\}$ as the set of non-negative integers, where $n\mathbb{N}_0$ denotes all elements in $\mathbb{N}_0$ multiplied by $n$. Additionally, $\mathbb{Z}$, $\mathbb{Q}$, $\mathbb{R}$, and $\mathbb{C}$ represent the sets of integers, rational, real, and complex numbers, respectively.\\\\
We denote the gamma \cite[(5.2.1)]{bib23}, digamma \cite[(5.2.2)]{bib23}, tetragamma, pentagamma, and hexagamma functions \cite[\S5.15]{bib23} of argument $z$ as $\Gamma(z)$, $\psi(z)$, $\psi_2(z)$, $\psi_3(z)$, and $\psi_4(z)$, respectively. Here, $\psi_n(z)$ is the polygamma function, defined as the $n$-th derivative of $\ln{\Gamma(z)}$, and $n\in\mathbb{N}_0$. The digamma function can be expressed as \cite[\S1.7(6)]{bib23}:
\[\psi(z) = -\gamma +\sum_{k=0}^\infty \left(\frac{1}{k+1} - \frac{1}{k+z}\right), \quad z \in \mathbb{C}\setminus- \mathbb{N}_0.\]
For positive integer values of $z$, the digamma function simplifies to \cite[\S1.7.1(9)]{bib2}:
\begin{equation}\label{harmonicn}
\psi(k+1) = -\gamma + H_k, \quad k \in \mathbb{N}.
\end{equation}
The recurrence relation for the digamma function is given by \cite[(5.5.2)]{bib23}:
\begin{equation}\label{recc}
\psi(z+1) = \psi(z) + \frac{1}{z}.
\end{equation}
The duplication formula for the $\psi(z)$ is \cite[(5.5.8)]{bib23}:
\begin{equation}\label{dupl}
\psi\left(z + \frac{1}{2}\right) = 2\psi(2z) - \psi(z) - \ln{4}, \quad z \in \C \setminus-\mathbb{N}_0.
\end{equation}
The Lerch transcendent is defined as \cite[(24.14.1)]{bib23}:
$$\Phi(z, s, a) = \sum_{n=0}^{\infty} \frac{z^n}{(n + a)^s}, \quad \lvert z\rvert \leq 1, \Re\,s > 1, a \not\in -\mathbb{N}_0.$$
The polygamma function can be expressed as $\psi_n(z) = (-1)^{n-1} n! \Phi(1, n+1, z)$, where $n \in \N$. The Dirichlet eta function is defined as $\eta(n) := \Phi(-1, n, 1)$, where $\Re n > 0$ \cite[\S1.12(2)]{bib23}. The Riemann zeta function \cite[\S 25.2]{bib23} and the Hurwitz zeta function \cite[\S 25.11]{bib23} are, respectively, defined as:
$$\zeta(s) = \sum_{n=1}^{\infty} \frac{1}{n^s}, \quad \zeta(s, z) = \sum_{n=0}^{\infty} \frac{1}{(n+z)^s},$$
where $z \not\in -\mathbb{N}_0$, $\Re\,s > 1$. The domain $\Re\,s > 1$ of the Riemann and Hurwitz zeta functions can be extended to $s \in \mathbb{C}\setminus \{1\}$ through analytic continuation, using for instance, the Hermite integral representation for the Hurwitz zeta function \cite[(25.11.29)]{bib23} 
\begin{equation}\label{herm}
\zeta(s, z) = \frac{z^{-s}}{2} + \frac{z^{1-s}}{s-1} + 2\int_0^\infty \frac{\sin\left(s\arctan\left(x/z\right)\right)}{\left(x^2 + z^2\right)^{\frac{s}{2}}\left(\mathrm{e}^{2\pi x} - 1\right)}\, \mathrm{d}x.
\end{equation}
The relationship between the Dirichlet eta function and the Riemann zeta function is given by $\eta(n) = \left(1- 2^{1-n}\right)\zeta(n)$ \cite[\S1.12(2)]{bib2}. The polylogarithm function \cite[\S 25.12(ii), \S 25.14.3]{bib23}, $\li_s(z)$, is defined as: $\li_s(z) = z\Phi(z, s, 1)$, where $\Re s > 1,\, |z| \leq 1$. The dilogarithm function, $\li_2(z)$, has the integral representation \cite[(25.12.2)]{bib23}
\begin{equation}\label{intrep}
\li_2(z) = -\int_0^z \frac{\ln(1-t)}{t}\, \rmd{d}t, \quad z \in \C \setminus(1, \infty).
\end{equation}
The domain $|z| \leq 1$ can be extended to the entire complex plane through analytic continuation. This can be achieved using, for instance, the integral representation \cite[(25.14.6)]{bib23}
\[\li_s(z) = \frac{1}{2}z + \int_0^\infty \frac{z^{t+1}}{(1+t)^s} \, \rmd{d}t - 2z\int_0^\infty \frac{\sin\left(t\ln{z} - s\arctan{t}\right)}{\left(1+t^2\right)^{\frac{s}{2}} \left(e^{2\pi t} - 1\right)} \, \rmd{d}t, \quad \Re(s) > 0\,\,\textup{if}\,\, z\in \C\setminus[1, \infty).\]
Throughout this work, we utilize the property that $\li_s(z)$ is defined for all complex $z$.
\section{Results}
In this section, we present the main findings and outcomes of our study. We begin by deriving a closed form for \eqref{ruh1} using known results. Afterwards, we present a novel generalization of this integral.

\begin{lemma}[V\u{a}lean] \label{lemma4}
The following results are valid:
\begin{align}
&\sum_{k=1}^\infty \frac{H_k}{k^3} = \frac{\pi^4}{72}, \label{hrmkcu}\\
&\sum_{k=1}^\infty \frac{H_k}{(2k-1)^3} =  -\frac{\pi^2}{4} + \frac{\pi^4}{64} + 2\ln{2} + \frac{7\zeta(3)}{4} - \frac{7\ln{2}\zeta(3)}{4}, \label{ftpv}\\
& \sum_{k=1}^\infty \frac{H_{2k-1}}{(2k-1)^3} = \frac{\pi^4}{45} + \frac{\pi^2}{24}\ln^2{2} - \frac{\ln^4{2}}{24} - \frac{7\ln{2}\zeta (3)}{8} - \li_4\left(\frac{1}{2}\right), \label{wgt1} \\
& \sum_{k=1}^\infty \frac{H_{2k}}{k^3} =  -\frac{\pi ^4}{15}-\frac{\pi ^2 }{3}\ln^2{2} +\frac{\ln ^4{2}}{3} +7 \ln{2} \zeta (3) + 8 \li_4\left(\frac{1}{2}\right).\label{relf1}
\end{align}
\end{lemma}

\begin{proof}
The first series \eqref{hrmkcu} follows from \cite[\S 6.19,  pp.~601, (6.149)]{bib9} for the case $p=3$. Using $H_k = H_{k+1} - \frac{1}{k+1}$ in \cite[\S 4.19,  pp.~420, (4.102)]{bib9} for the case $p=3$ and reindexing, the proof of \eqref{ftpv} is complete. Adding \eqref{hrmkcu} to the closed form of $\sum_{k=1}^\infty \frac{(-1)^{k-1}H_{k}}{k^3}$ provided in \cite[\S 6.52,  pp.~502]{bib32}, the proof of \eqref{wgt1}. Subtracting both series, the closed form of \eqref{relf1} is complete.
\end{proof}

\begin{remark}
The closed form for \eqref{ruh1} is
\begin{equation}\label{qot1}
\int_0^1 \frac{\li_2(z) \ln(1+z)}{z} \, \rmd{d}z = \frac{\ln^4{2}}{12} - \frac{\pi^2 \ln^2{2}}{12} - \frac{\pi^4}{60} + \frac{7\ln{2}\zeta(3)}{4} + 2\li_4\left(\frac{1}{2}\right).
\end{equation}
\end{remark}

\begin{proof}
By utilizing the series representation of $\li_2(z)$ for $|z| \leq 1$, we can express \eqref{ruh1} as
\begin{equation}\label{abv2}
\int_0^1 \frac{\li_2(z) \ln(1+z)}{z} \, \rmd{d}z = \sum_{k=1}^\infty \frac{1}{k^2} \int_0^1 z^{k-1} \ln(1+z) \, \rmd{d}z.
\end{equation}
Upon integrating term by term, we obtain
\begin{equation}\label{abv1}
\int_0^1 z^{k-1} \ln(1+z) \, \rmd{d}z = \frac{\ln{2}}{k} - \frac{1}{k}\int_0^1 \frac{z^k}{1+z} \, \rmd{d}z =  \frac{\ln{2}}{k} - \frac{1}{2k}\left(\psi\left(\frac{k+2}{2}\right) - \psi\left(\frac{k+1}{2}\right)\right).
\end{equation}
Applying \eqref{dupl} and subsequently employing \eqref{harmonicn} in \eqref{abv1}, we derive
\begin{equation}\label{abv3}
\int_0^1 z^{k-1} \ln(1+z) \, \rmd{d}z  = \frac{H_k - H_{\frac{k}{2}}}{k}.
\end{equation}
Substituting \eqref{abv3} into \eqref{abv2} and applying \eqref{harmonicn} and \eqref{dupl} after splitting into odd and even parts, we obtain
\begin{equation}\label{hafe2}
\begin{split}
\int_0^1 \frac{\li_2(z) \ln(1+z)}{z} \, \rmd{d}z &= \sum_{k=1}^\infty  \frac{H_k - H_{\frac{k}{2}}}{k^3} 
\\&= -2 + \frac{7}{8} \sum_{k=1}^\infty \frac{H_k}{k^3} + \frac{7\ln{2}\zeta(3)}{4}  - \sum_{k=1}^\infty \frac{2H_{2k+1} - H_k}{(2k+1)^3}.
\end{split}
\end{equation}
Reindexing the series \eqref{ftpv}, we have
\begin{equation}\label{hafe1}
\sum_{k=1}^\infty \frac{H_k}{(2k+1)^3} = \frac{\pi^4}{64} - \frac{7\ln{2}}{4} \zeta(3).
\end{equation}
Substituting \eqref{hrmkcu}, \eqref{wgt1}, and \eqref{hafe1} into \eqref{hafe2}, we successfully conclude the proof of \eqref{qot1}. Alternatively, we can employ the series representation of $\ln(1 + z)$ for $|z| \leq 1$ to establish the proof of \eqref{qot1}. This approach yields
\begin{equation}\label{haffn}
\begin{split}
\int_0^1 \frac{\li_2(z) \ln(1+z)}{z} \, \rmd{d}z &= \sum_{k=1}^\infty \frac{(-1)^{k-1}}{k}\left(\frac{\pi^2}{6k} - \frac{H_k}{k^2}\right) = \frac{\pi^4}{72} - \sum_{k=1}^\infty \frac{(-1)^{k-1}H_k}{k^3} 
\\&= \sum_{k=1}^\infty \frac{H_k}{k^3} - \frac{(-1)^{k-1}H_k}{k^3} = \frac{1}{4}\sum_{k=1}^\infty \frac{H_{2k}}{k^3}.
\end{split}
\end{equation}
By substituting \eqref{relf1} into \eqref{haffn}, we readily conclude the proof of \eqref{qot1}. Notably, this approach appears to be faster, as we immediately recognise $\frac{\pi^4}{72}$ as the closed form of $\sum_{k=1}^\infty \frac{H_k}{k^3}$.
\end{proof}

\subsection{Generalization of the dilogarithmic integral and related integrals} \label{sec3.1}
In the following theorems, we provide the generalization \eqref{begeqh1}, and the generalization of integrals related to \eqref{begeqh1}. The closed forms for integrals presented in Theorems \ref{thmabd1}--\ref{thmhaf5} are new and have not been presented elsewhere in the literature. In these theorems, we avoid computations of logarithm of negative real numbers.
\begin{theorem} Let $a \in \C \setminus(-\infty, -1)$. Then \label{thmabd1}
\begin{equation}\label{thmhaf1}
\begin{split}
\int_0^1 \frac{\ln{z} \ln(1 + az) \ln(1 - z)}{z} \, \rmd{d}z + \frac{1}{2} \int_0^1 \frac{\ln{z}\ln^2(1 + a z)}{z}&= - \frac{\left(\li_2(-a)\right)^2}{2}+ \frac{\pi^2}{6}\li_2(-a)  
\\&\quad-2\li_4(-a).
\end{split}
\end{equation}
\begin{equation}\label{thmhaf2}
\int_0^1 \frac{\li_2(z) \ln(1 + az)}{z}\, \rmd{d}z + \frac{1}{2}\int_0^1 \frac{\ln{z} \ln^2(1 + az)}{z} \, \rmd{d}z = - \frac{\pi^2}{6}\li_2(-a) + \li_4(-a) .
\end{equation}
\end{theorem}

\begin{proof}
In an effort to circumvent the computation of logarithm of negative real numbers, the initial part of the proof addresses the case where $a \in \C \setminus (-\infty, 0)$, while the subsequent section pertains to the scenario where $a \in \C \setminus (-\infty, -1) \cup (0, \infty)$. In both cases, it follows that $a \in \C \setminus(-\infty, -1)$. Now, we begin by performing integration by parts, resulting in
\begin{equation}\label{major1}
\int_0^1 \frac{\li_2(z) \ln(1 + az)}{z}\, \rmd{d}z = \int_0^1 \frac{\ln{z} \ln(1 + az) \ln(1 - z)}{z} \, \rmd{d}z - a\int_0^1 \frac{\ln{z} \li_2(z)}{1 + az} \, \rmd{d}z.
\end{equation}
V\u{a}lean employed the Cauchy product of two series \cite{bibch} to derive \cite[\S 3.49, pp.~335, (3.333)]{bib9}
\begin{equation}\label{cauchy1}
\left(\li_2(a)\right)^2 = 4\sum_{k=1}^\infty \frac{a^k H_k}{k^3} + 2\sum_{k=1}^\infty \frac{a^k H_k^{(2)}}{k^2} - 6 \sum_{k=1}^\infty \frac{a^k}{k^4}, \quad a \in \C, |a| \leq 1,
\end{equation}
Utilizing \eqref{cauchy1}, V\u{a}lean \cite[\S1.49, (1.218)]{bib9} provided the closed form for the second resulting integral in \eqref{major1}, with $a$ in \eqref{major1} substituted with $-a$. This yields
\begin{equation}\label{major2}
\int_0^1 \frac{\ln{z} \li_2(z)}{1 + az} \, \rmd{d}z = -\frac{\left(\li_2(-a)\right)^2}{2a}  + \frac{\pi^2\li_2(-a)}{3a} - \frac{3\li_4(-a)}{a}, \quad a \in \C\setminus(-\infty, -1) \cup \{0\}.
\end{equation}
At this point, our focus narrows down to obtaining an expression for the first resulting integral in \eqref{major1}. By carrying out term-by-term integration, we deduce
\begin{equation}\label{major3}
\begin{split}
\int_0^1 \frac{\ln{z} \ln(1 + az) \ln(1 - z)}{z} \, \rmd{d}z &= \sum_{k=1}^\infty \frac{(-1)^{k-1} a^k }{k} \int_0^1 z^{k-1} \ln{z} \ln(1 - z) \, \rmd{d}z
\\&= -\li_4(-a) - \sum_{k=1}^\infty \frac{(-1)^k a^k H_k}{k^3} + \sum_{k=1}^\infty \frac{(-1)^k a^k \psi_1(k)}{k^2}.
\end{split}
\end{equation}
By implementing the notation change $H^{(2)}_n = \pi^2/ 6- \psi_1(n+1)$ in \eqref{cauchy1} and rearranging, we have
\begin{equation}\label{cauchy2}
\sum_{k=1}^\infty \frac{(-1)^k a^k \psi_1(k)}{k^2} =  \frac{\pi^2}{6}\li_2(-a) - 2\li_4(-a) - \frac{1}{2}\left(\li_2(-a)\right)^2 + 2\sum_{k=1}^\infty \frac{(-1)^k a^k H_k}{k^3}.
\end{equation}
Upon substituting \eqref{cauchy2} into \eqref{major3}, we derive
\begin{equation}\label{major4}
\begin{split}
\int_0^1 \frac{\ln{z} \ln(1 + az) \ln(1 - z)}{z} \, \rmd{d}z &= -3\li_4(-a) + \frac{\pi^2}{6}\li_2(-a)  - \frac{1}{2}\left(\li_2(-a)\right)^2 
\\&\qquad+ \sum_{k=1}^\infty \frac{(-1)^k a^k H_k}{k^3}.
\end{split}
\end{equation}
We employ the relationship \cite[\S4.6, pp.~504]{bib9}
\begin{equation}\label{major5}
\sum_{k=1}^\infty \frac{H_k}{k+1} p^{k} = \frac{\ln^2(1 - p)}{2p}, \quad p \in \C, |p| \leq 1, p \neq 0, 1.
\end{equation}
By substituting $p$ with $-p$ in \eqref{major5}, integrating over the interval from $p = 0$ to $z$, we deduce
\begin{equation}\label{major6}
\sum_{k=1}^\infty \frac{(-1)^{k+1}H_k}{(k+1)^2} z^{k+1} = \frac{1}{2} \int_0^z \frac{\ln^2(1 + p)}{p}\, \rmd{d}p, \quad z \in \C \setminus (-\infty, 0) \cup (1, \infty).
\end{equation}
Subsequently, we perform a change of variable from $z$ to $t$ in \eqref{major6}, divide by $t$, and integrate once more from $t = 0$ to $a$, resulting in
\begin{equation}\label{major7}
\sum_{k=1}^\infty \frac{(-1)^{k+1}H_k}{(k+1)^3} a^{k+1} = \frac{1}{2} \int_0^a \int_0^t \frac{\ln^2(1 + p)}{p}\, \rmd{d}p \, \rmd{d}t.
\end{equation}
Substituting \cite[\S 1.4, pp.~3, (1.9)]{bib32} into \eqref{major7} and reindexing the series on the left-hand side, we arrive at
\begin{equation}\label{major8}
\begin{split}
\sum_{k=1}^\infty \frac{(-1)^{k-1} H_k a^k}{k^3}&= \int_0^a  \frac{1}{t}\left(-\li_3(-t) - \zeta(3) + \frac{\ln^3(1 + t)}{3} + \li_3\left(\frac{1}{1+t}\right) \right.
\\&\qquad \left.+ \ln(1+ t)\li_2\left(\frac{1}{1+t}\right)  - \frac{1}{2}\ln{t}\ln^2(1 + t)\right) \rmd{d}t, \quad a \in \C \setminus (-\infty, 0) \cup (1, \infty).
\end{split}
\end{equation}
Through integration by parts, we obtain these two generalized results:
\begin{equation}\label{combh1}
\begin{split}
\int_0^a  \frac{\zeta(3) -  \li_3\left(\frac{1}{1+t}\right)}{t} \, \rmd{d}t &= \ln{a}\left(\zeta(3) -  \li_3\left(\frac{1}{1+a}\right)\right)  - \ln{a}\ln(1 + a) \li_2\left(\frac{1}{1+a}\right)
\\&\quad- \li_2\left(\frac{1}{1+a}\right) \li_2(-a) + \ln{a}\ln(1 + a) \li_2(-a)   
\\&\quad- \frac{1}{2}\ln^2(1+ a)\li_2(-a)   + \frac{1}{2}\left(\li_2(-a)\right)^2 + \frac{\ln^2{a}}{2}\ln(1 + a) 
\\&\quad- \frac{\ln{a}}{3}\ln^3(1 + a) - \frac{1}{6} \int_0^a \frac{\ln^3(1 + t)}{t} \, \rmd{d}t,
\end{split}
\end{equation}
\begin{equation}\label{combh2}
\begin{split}
\int_0^a  \frac{\ln(1 + t) \li_2\left(\frac{1}{1+t}\right)}{t} \, \rmd{d}t &= -\li_2\left(\frac{1}{1+a}\right) \li_2(-a) + \li_2(-a)\ln{a}\ln(1 + a) 
\\&\quad + \frac{1}{2}\left(\li_2(-a)\right)^2 - \frac{\ln^2(1  + a)}{2} \li_2(-a) 
\\&\quad- \frac{1}{2} \int_0^a \frac{\ln^3(1 + t)}{t} \, \rmd{d}t + \int_0^a \frac{\ln{t}\ln^2(1 + t)}{t}\, \rmd{d}t.
\end{split}
\end{equation}
By employing \eqref{combh1} and \eqref{combh2} within \eqref{major8}, we obtain
\begin{equation}\label{major9}
\begin{split}
\sum_{k=1}^\infty \frac{(-1)^{k-1} H_k a^k}{k^3} &= -\li_4(-a) + \ln{a}\left(\li_3\left(\frac{1}{1+a}\right) - \zeta(3)\right) + \ln{a}\ln(1 + a) 
\\&\quad \times \li_2\left(\frac{1}{1+a}\right) + \ln{a}\ln(1 + a)\left(\frac{\ln^2(1 + a)}{3} - \frac{\ln{a}}{2}\right) 
\\&\quad+ \frac{1}{2}\int_0^a \frac{\ln{t}\ln^2(1 + t)}{t}\, \rmd{d}t.
\end{split}
\end{equation}
Utilizing \cite[\S 1.4, pp.~3, (1.9)]{bib32} in the last integral in \eqref{major9}, we arrive at
\begin{equation}\label{major10}
\begin{split}
\frac{1}{2}\int_0^a \frac{\ln{t}\ln^2(1 + t)}{t}\, \rmd{d}t &= \ln{a}\left(\zeta(3) - \li_3\left(\frac{1}{1+a}\right)\right) + \frac{1}{2}\ln^2{a}\ln(1 + a)
\\&\quad - \frac{1}{3}\ln{a}\ln^3(1 + a) - \ln{a}\ln(1 + a)\li_2\left(\frac{1}{1+a}\right) 
\\&\quad+ \frac{1}{2} \int_0^1 \frac{\ln{t}\ln^2(1 + a t)}{t}\, \rmd{d}t.
\end{split}
\end{equation}
Substituting \eqref{major10} into \eqref{major9}, we deduce
\begin{equation}\label{major11}
\sum_{k=1}^\infty \frac{(-1)^{k} H_k a^k}{k^3} = \li_4(-a)  - \frac{1}{2} \int_0^1 \frac{\ln{z}\ln^2(1 + a z)}{z}\, \rmd{d}z, \quad a \in \C \setminus (-\infty, 0) \cup (1, \infty).
\end{equation}
Likewise, for $a \in \C \setminus (-\infty, 0] \cup [1, \infty)$, we obtain using \eqref{major5}
\begin{equation}\label{simha1}
\int_0^a \frac{\ln{z}\ln^2(1 - z)}{z}\, \rmd{d}z = - 2\ln{a}\li_3(a)  + 2\li_4(a)  -2\sum_{k=1}^\infty \frac{H_k a^k}{k^3}  + 2\ln{a}\sum_{k=1}^\infty \frac{H_k a^k}{k^2}.
\end{equation}
By making use of  \cite[\S4.6, pp.~399, (4.36)]{bib9}
\begin{equation}\label{harmo1}
\sum_{k=1}^\infty \frac{H_k}{k^2} z^{k} = \zeta(3) + \li_2(1-z)\ln(1-z) + \li_3(z) - \li_3(1-z) + \frac{1}{2}\ln{z}\ln^2(1-z),
\end{equation}
in \eqref{simha1}, we deduce
\begin{equation}\label{simha2}
\sum_{k=1}^\infty \frac{H_k a^k}{k^3} = \li_4(a) -  \frac{1}{2} \int_0^1 \frac{\ln{z}\ln^2(1 - a z)}{z}\, \rmd{d}z, \quad a \in \C \setminus (-\infty, 0) \cup (1, \infty).
\end{equation}
Substituting $a$ with $-a$ in \eqref{simha2}, we once again arrive at \eqref{major11}. However, this time, \eqref{major11} is applicable for $a \in \C \setminus (-\infty, -1) \cup (0, \infty)$. Further substituting \eqref{major11} into \eqref{major4}, we derive
\begin{equation}\label{major12}
\begin{split}
\int_0^1 \frac{\ln{z} \ln(1 + az) \ln(1 - z)}{z} \, \rmd{d}z &= -2\li_4(-a) + \frac{\pi^2}{6}\li_2(-a)  - \frac{1}{2}\left(\li_2(-a)\right)^2 
\\&\quad-\frac{1}{2} \int_0^1 \frac{\ln{z}\ln^2(1 + a z)}{z}\, \rmd{d}z.
\end{split}
\end{equation}
Rearranging \eqref{major12} concludes the proof of \eqref{thmhaf1}. Finally, by substituting \eqref{major2} and \eqref{major12} into \eqref{major1}, we complete the proof of \eqref{thmhaf2}.
\end{proof}

\begin{remark}
In Theorem \ref{thmabd1}, we have the following special values:
\begin{equation}
\int_0^1 \frac{\ln{z} \ln(1 + z) \ln(1 - z)}{z} \, \rmd{d}z + \frac{1}{2} \int_0^1 \frac{\ln{z}\ln^2(1 + z)}{z} = \frac{\pi^4}{480},
\end{equation}
\begin{equation}
\int_0^1 \frac{\li_2(z) \ln(1 + z)}{z}\, \rmd{d}z + \frac{1}{2}\int_0^1 \frac{\ln{z} \ln^2(1 + z)}{z} \, \rmd{d}z = \frac{\pi^4}{240},
\end{equation}
\begin{equation}
\int_0^1 \frac{\li_2(z) \ln(1 - z)}{z}\, \rmd{d}z + \frac{1}{2}\int_0^1 \frac{\ln{z} \ln^2(1 - z)}{z} \, \rmd{d}z = -\frac{\pi^4}{60}.
\end{equation}
\end{remark}

\begin{theorem}Let $a \in \C \setminus(-\infty, 0]$. Then \label{thmhaf22}
\begin{equation*} 
\begin{split}
\int_0^1 \frac{\ln{z} \ln(1 + az) \ln(1 - z)}{z} \, \rmd{d}z &= -\frac{\pi^4}{90} - \frac{\left(\li_2(-a)\right)^2}{2}+ \frac{\pi^2}{6}\li_2(-a) - \li_4(-a)  
\\&\quad+ \frac{\pi^2}{12}  \ln^2(1 + a) + \frac{1}{3}\ln{a}\ln^3(1 + a) - \frac{1}{4}\ln^4(1 + a) 
\\&\quad + \ln(1 + a)\left(\li_3\left(\frac{1}{1+a}\right) + \li_3\left(\frac{a}{1+a}\right)\right) 
\\&\quad+ \li_4\left(\frac{1}{1+a}\right) + \li_4\left(\frac{a}{1+a}\right).
\end{split}
\end{equation*}
\end{theorem}

\begin{proof}
We establish the proof by evaluating the second integral in \eqref{thmhaf1}. By algebraic substitutions, we have
\begin{equation}\label{disin}
\int_0^1 \frac{\ln{z}\ln^2(1 + az)}{z}\, \rmd{d}z = \int_{\frac{1}{1+a}}^1 \frac{\ln(1 - z)\ln^2{z} - \ln^3{z} - \ln{a}\ln^2{z}}{z(1 - z)} \, \rmd{d}z.
\end{equation}
For the first resulting integral, we have
\begin{equation} \label{eq365h}
\begin{split}
\int_{\frac{1}{1+a}}^1 \frac{\ln(1 - z)\ln^2{z}}{z(1- z)} &= \int_0^1 \frac{\ln^2{z}\ln(1 - z)}{z} \, \rmd{d}z - \int_0^{\frac{1}{1+a}} \frac{\ln^2{z}\ln(1 - z)}{z} \, \rmd{d}z 
\\&\quad+  \int_0^{\frac{a}{1+a}} \frac{\ln{z}\ln^2(1 - z)}{z} \, \rmd{d}z.
\end{split}
\end{equation}
Using \eqref{harmo1} and \cite[\S4.6, pp.~399, (4.38)]{bib9}
\begin{equation}\label{harmo2}
\begin{split}
\sum_{k=1}^\infty \frac{H_k}{k^3} z^{k} &=  \frac{\pi^4}{90} + \zeta(3)\ln(1 - z) + \frac{\pi^2}{12}\ln^2(1 - z) + \frac{1}{24}\ln^4(1 - z) - \frac{1}{6}\ln{z}\ln^3(1 - z) 
\\&\quad- \ln(1 - z)\li_3(z) + 2\li_4(z) - \li_4(1-z) + \li_4\left(\frac{z}{z-1}\right)
\end{split}
\end{equation}
in \eqref{simha1}, we deduce
\begin{equation}\label{newinh1}
\begin{split}
\int_0^z \frac{\ln{t}\ln^2(1 - t)}{t} \, \rmd{d}t &= -\frac{\pi^4}{45} - 2\zeta(3)\ln(1 - z) - \frac{\pi^2}{6}\ln^2(1 - z) - \frac{1}{12}\ln^4(1 - z) 
\\&\quad+ \frac{1}{3}\ln{z}\ln^3(1 - z) +  2\ln(1 - z)\li_3(z) - 2\li_4(z) + 2\li_4(1 - z)
\\&\quad- 2\li_4\left(\frac{z}{z-1}\right) + 2\ln{z}\zeta(3) + 2\ln{z}\ln(1 - z)\li_2(1 - z) 
\\&\quad- 2\ln{z}\li_3(1 - z) + \ln^2{z}\ln^2(1 - z).
\end{split}
\end{equation}
The following results are obtained through straightforward term-by-term integration
\begin{equation}\label{frsin1}
\begin{split}
\int_0^{\frac{1}{1+a}} \frac{\ln^3{z}}{1 - z} \, \rmd{d}z &= (\ln{a} - \ln(1+ a)) \ln^3(1 + a)- 3\ln^2(1 + a)\li_2\left(\frac{1}{1+a}\right) 
\\&\quad- 6\ln(1 + a)\li_3\left(\frac{1}{1+a}\right)  - 6\li_4\left(\frac{1}{1+a}\right),
\end{split}
\end{equation}
\begin{equation}\label{frsin2}
\begin{split}
\int_0^{\frac{a}{1+a}} \frac{\ln^2(1 - z)}{z}\, \rmd{d}z &= 2\zeta(3) - 2\li_3\left(\frac{1}{1+a}\right) + (\ln{a} - \ln(1 +a))\ln^2(1 + a) 
\\&\quad- 2\ln(1+ a)\li_2\left(\frac{1}{1+a}\right).
\end{split}
\end{equation}
By replacing $z$ with $a/(a+1)$ in \eqref{newinh1} and employing and \eqref{frsin1}, and further substituting the latter in \eqref{disin} while taking into account \eqref{frsin2}, we ultimately arrive at
\begin{equation}\label{fimp1}
\begin{split}
\int_0^1 \frac{\ln{z}\ln^2(1 + az)}{z}\, \rmd{d}z &= \frac{\pi^4}{45} - \frac{\pi^2}{6} \ln^2(1 + a) - \frac{2}{3}\ln{a}\ln^3(1 + a) + \frac{1}{2}\ln^4(1 + a) 
\\&\quad- 2\ln(1 + a)\left(\li_3\left(\frac{1}{1+a}\right) + \li_3\left(\frac{a}{1+a}\right)\right) - 2\li_4(-a)
\\&\quad - 2\li_4\left(\frac{1}{1+a}\right) - 2\li_4\left(\frac{a}{1+a}\right), \quad a \in \C \setminus (-\infty, 0].
\end{split}
\end{equation}
Substituting \eqref{fimp1} into \eqref{thmhaf1}, we thereby conclude the proof of Theorem \ref{thmhaf22}.
\end{proof}

\begin{theorem} Let $a \in \C \setminus(-\infty, 0]$. Then \label{thmhaf3}
\begin{equation} \label{inthmhaf3}
\begin{split}
\int_0^1 \frac{\li_2(z) \ln(1 + az)}{z}\, \rmd{d}z  &= -\frac{\pi^4}{90} - \frac{\pi^2}{6}\li_2(-a) + 2\li_4(-a) +  \frac{\pi^2}{12}  \ln^2(1 + a) 
\\&\quad+  \frac{1}{3}\ln{a}\ln^3(1 + a) - \frac{1}{4}\ln^4(1 + a) + \ln(1 + a)\li_3\left(\frac{1}{1+a}\right) 
\\&\quad+ \ln(1 + a)\li_3\left(\frac{a}{1+a}\right)+ \li_4\left(\frac{1}{1+a}\right) + \li_4\left(\frac{a}{1+a}\right).
\end{split}
\end{equation}
\end{theorem}

\begin{proof}
Through the substitution of \eqref{fimp1} into \eqref{thmhaf2}, we establish the proof of Theorem \ref{thmhaf3}.
\end{proof}

\begin{remark}Considering \eqref{major4}, we substitute $z=-a$ in \eqref{harmo2}, resulting in \label{greatremark}
\begin{equation}\label{harmo21}
\begin{split}
\sum_{k=1}^\infty \frac{(-1)^{k} H_k a^k}{k^3} &=  \frac{\pi^4}{90} + \zeta(3)\ln(1 + a) + \frac{\pi^2}{12}\ln^2(1 + a) + \frac{1}{24}\ln^4(1 + a) 
\\&\quad - \frac{1}{6}\ln{(-a)}\ln^3(1 + a)  - \ln(1 + a)\li_3(-a) + 2\li_4(-a) - \li_4(1+a) 
\\&\quad + \li_4\left(\frac{a}{a+1}\right).
\end{split}
\end{equation}
While V\u{a}lean impressively provided different representations for \eqref{harmo1} and \eqref{harmo2}, it is worth noting that the logarithms on the right-hand side of  \eqref{harmo1} and \eqref{harmo2} avoid negative values for $z \in \C\setminus(-\infty, 0] \cup [1, \infty)$. Consequently, \eqref{harmo21} avoids logarithms of negative numbers for $a \in \C\setminus (-\infty, -1] \cup [0, \infty)$. In the other region where $a \in \C\setminus(-\infty, 0] \cup [1, \infty)$, by employing the relationships \eqref{major11} and \eqref{fimp1} that we have established, we arrive at
\begin{equation}\label{news1}
\begin{split}
\sum_{k=1}^\infty \frac{(-1)^{k} H_k a^k}{k^3} &= -\frac{\pi^4}{90} + \frac{\pi^2}{12} \ln^2(1 + a) + \frac{1}{3}\ln{a}\ln^3(1 + a) + \frac{1}{2}\ln^4(1 + a) 
\\&\quad +\ln(1 + a)\left(\li_3\left(\frac{1}{1+a}\right) + \li_3\left(\frac{a}{1+a}\right)\right) + 2\li_4(-a)
\\&\quad + \li_4\left(\frac{1}{1+a}\right) + \li_4\left(\frac{a}{1+a}\right),
\end{split}
\end{equation}
which is valid for all $z \in \C$, where $|z| \leq 1$ and $z \neq -1, 0$. Similarly, from \eqref{major6} and \cite[\S 1.4, pp.~3, (1.9)]{bib32}, we derive
\begin{equation} \label{lastheq1}
\begin{split}
\sum_{k=1}^\infty \frac{(-1)^k H_k}{(k+1)^2} z^{k+1} &= -\zeta(3) + \frac{1}{3}\ln^3(1 + z) + \li_3\left(\frac{1}{1+z}\right) + \ln(1+z)\li_2\left(\frac{1}{1+z}\right)
\\&\quad-\frac{1}{2}\ln{z}\ln^2(1+z).
\end{split}
\end{equation}
Reindexing the series in the left-hand side of \eqref{lastheq1}, we have
\begin{equation} \label{lastheq2}
\sum_{k=1}^\infty \frac{(-1)^k H_k}{(k+1)^2} z^{k+1} = -\sum_{k=1}^\infty \frac{(-1)^k H_{k}}{k^2} z^{k} + \li_3(-z).
\end{equation}
Substituting \eqref{lastheq2} into \eqref{lastheq1}, we obtain
\newpage
\begin{equation} \label{stah12}
\begin{split}
\sum_{k=1}^\infty \frac{(-1)^k H_{k}}{k^2} z^{k} &= \zeta(3)  - \frac{1}{3}\ln^3(1 + z) + \li_3(-z) - \li_3\left(\frac{1}{1+z}\right) - \ln(1+z)\li_2\left(\frac{1}{1+z}\right)
\\&\quad+\frac{1}{2}\ln{z}\ln^2(1+z).
\end{split}
\end{equation}

The representations \eqref{news1} and \eqref{stah12}, applicable to $a \in \C\setminus(-\infty, 0] \cup [1, \infty)$, circumvents the need to compute logarithms of negative numbers. The advantage of \eqref{news1} over \eqref{harmo21} lies in its application to the evaluation of the integral in \eqref{major4}. Specifically, when considering avoidance of logarithms of negative numbers and the analytic continuation of the polylogarithm function, \eqref{news1} allows for a more general result, extending the domain of the integral from $\C\setminus(-\infty, 0] \cup [1, \infty)$ to $\C\setminus(-\infty, 0]$. Conversely, \eqref{harmo21} still restricts the domain of the integral to $\C\setminus (-\infty, -1] \cup [0, \infty)$, avoiding logarithms of negative numbers while still accommodating the analytic continuation of the polylogarithm function. In light of these observations, we reveal two additional closed forms for the integrals in Theorems \ref{thmhaf22} and \ref{thmhaf3} in the following theorems.
\end{remark}

\begin{theorem} Let $a \in \C\setminus (-\infty, -1] \cup [0,  \infty)$. Then \label{thmhaf4}
\begin{equation} 
\begin{split}
\int_0^1 \frac{\ln{z} \ln(1 + az) \ln(1 - z)}{z} \, \rmd{d}z &=  \frac{\pi^4}{90} -\li_4(-a) + \frac{\pi^2}{6}\li_2(-a)  - \frac{1}{2}\left(\li_2(-a)\right)^2 
\\&\quad+ \zeta(3)\ln(1 + a) + \frac{\pi^2}{12}\ln^2(1 + a) + \frac{1}{24}\ln^4(1 + a) 
\\&\quad- \frac{1}{6}\ln{(-a)}\ln^3(1 + a) - \ln(1+ a)\li_3(-a)  
\\&\quad- \li_4(1+a)+ \li_4\left(\frac{a}{a+1}\right).
\end{split}
\end{equation}
\end{theorem}
\begin{proof}
Upon substituting \eqref{harmo21} into \eqref{major4}, we conclude the proof of Theorem \ref{thmhaf4}.
\end{proof}

\begin{theorem} Let $a \in \C\setminus (-\infty, -1] \cup [0,  \infty)$. Then \label{thmhaf5}
\begin{equation}
\begin{split}
\int_0^1 \frac{\li_2(z) \ln(1 + az)}{z}\, \rmd{d}z  &=  \frac{\pi^4}{90} + 2\li_4(-a) - \frac{\pi^2}{6}\li_2(-a)  + \zeta(3)\ln(1 + a) + \frac{\pi^2}{12}\ln^2(1 + a) 
\\&\quad + \frac{1}{24}\ln^4(1 + a) - \frac{1}{6}\ln{(-a)}\ln^3(1 + a) - \ln(1 + a)\li_3(-a) 
\\&\quad - \li_4(1+a) + \li_4\left(\frac{a}{a+1}\right).
\end{split}
\end{equation}
\end{theorem}

\begin{proof}
Considering \eqref{simha2} and \eqref{harmo21}, we infer the following closed form
\begin{equation}\label{thmhaf5eq}
\begin{split}
\int_0^1 \frac{\ln{z}\ln^2(1 + a z)}{z}\, \rmd{d}z &=  -\frac{\pi^4}{45} - 2\zeta(3)\ln(1 + a) - \frac{\pi^2}{6}\ln^2(1 + a) - \frac{1}{12}\ln^4(1 + a) 
\\&\quad + \frac{1}{3}\ln{(-a)}\ln^3(1 + a) +  2\ln(1 + a)\li_3(-a) - 2\li_4(-a) 
\\&\quad + 2\li_4(1+a) - 2\li_4\left(\frac{a}{a+1}\right), \quad a \in \C\setminus (-\infty, -1] \cup [0,  \infty).
\end{split}
\end{equation}
By substituting \eqref{thmhaf5eq} into \eqref{thmhaf2}, we conclude the proof of Theorem \ref{thmhaf5}.
\end{proof}

\noindent
In what follows, we apply our established identities to prove of Jonqui\`ere's inversion formula for specific cases. 
\begin{theorem}[Jonqui\`ere] Let $a \in \C\setminus (-\infty, -1] \cup [0,  \infty)$. Then \label{anathm}
\begin{equation}\label{myiden1}
\begin{split}
\li_4\left(\frac{z-1}{z}\right) + \li_4\left(\frac{z}{z-1}\right) &= -\frac{7\pi^4}{360} - \frac{1}{4}\ln^2{z}\ln^2(1 - z) - \frac{\pi^2}{12}\ln^2(1 - z) - \frac{1}{24}\ln^4(1 - z) 
\\&\quad+\frac{\pi^2}{6}\ln{z}\ln(1 - z) + \frac{1}{6}\ln{z}\ln^3(1 - z) - \frac{\pi^2}{12}\ln^2{z} 
\\&\quad + \frac{1}{6}\ln^3{z}\ln(1 - z) - \frac{1}{24}\ln^4{z}.
\end{split}
\end{equation}
\end{theorem}

\begin{proof}
On the transformation $z \longrightarrow 1 - z$, we have
\begin{equation}\label{eulan1}
\int_0^z \frac{\ln{t}\ln^2(1-t)}{t} \, \rmd{d}t + \int_0^{1-z} \frac{\ln{t}\ln^2(1-t)}{t} \, \rmd{d}t  = \frac{1}{2}\ln^2{z}\ln^2(1 - z) + \int_0^1 \frac{\ln{t}\ln^2(1-t)}{t} \, \rmd{d}t.
\end{equation}
This yields
\begin{equation}\label{eulan2}
\int_0^z \frac{\ln{t}\ln^2(1-t)}{t} \, \rmd{d}t + \int_0^{1-z} \frac{\ln{t}\ln^2(1-t)}{t} \, \rmd{d}t  = -\frac{\pi^4}{180} + \frac{1}{2}\ln^2{z}\ln^2(1 - z).
\end{equation}
Utilizing \eqref{newinh1} on the left-hand side of \eqref{eulan2}, we conclude the proof of Theorem \ref{anathm}.
\end{proof}

\begin{remark}
Let us define a function $F(z)$ as 
\[F(z) = \int_0^z \frac{\ln{t}\ln^2(1-t)}{t} \, \rmd{d}t.\]
With this notation, \eqref{eulan2} can be expressed as
\begin{equation}\label{eulan3}
F(z) + F(1 - z)  = -\frac{\pi^4}{180} + \frac{1}{2}\ln^2{z}\ln^2(1 - z).
\end{equation}
The identity \eqref{myiden1} is not new; however, we can refer to Theorem \ref{anathm} as a rediscovery, as the proof we provide is novel. Theorem \ref{anathm} can be derived by substituting 4 for $m$ and $\frac{z}{z-1}$ for $z$ in Jonqui\`ere's inversion formula \cite[\S1.11.1, pp.~31, (16)]{bib2}. Jonqui\`ere's inversion formula is derived from the Lerch transformation formula \cite[\S1.11(7), pp.~29]{bib2}, and the Lerch transformation formula is derived using the residue theorem (see \cite[\S1.11, pp.~28]{bib2}). This demonstrates the uniqueness of our proof.
\end{remark}

\begin{theorem}[Jonqui\`ere] \label{newthgh}
Let $z \in \C \setminus (-\infty, 0]$. Then
\begin{equation}
\li_4\left(-\frac{1}{a}\right) + \li_4(-a) =  -\frac{7\pi^4}{360} - \frac{\pi^2}{12}\ln^2{a} - \frac{\ln^4{a}}{24}.
\end{equation}
\end{theorem}

\begin{proof}
Define a function $G(z)$ as 
\[G(a) = \int_0^1 \frac{\ln{z}\ln^2(1 + az)}{z}\, \rmd{d}z.\]
Now, using \eqref{disin} and \eqref{eq365h}, we can establish the following relationship
\begin{equation}\label{geezh}
\begin{split}
G(a) + G\left(\frac{1}{a}\right) &=  -\frac{2\pi^4}{45} - \int_0^{\frac{1}{1+a}} \frac{\ln^2{z}\ln(1 - z)}{z} \, \rmd{d}z - \int_0^{\frac{a}{1+a}} \frac{\ln^2{z}\ln(1 - z)}{z} \, \rmd{d}z 
\\&\quad+  \int_0^{\frac{a}{1+a}} \frac{\ln{z}\ln^2(1 - z)}{z} \, \rmd{d}z +  \int_0^{\frac{1}{1+a}} \frac{\ln{z}\ln^2(1 - z)}{z} \, \rmd{d}z 
\\&\quad -\int_{\frac{1}{1+a}}^1 \frac{\ln^3{z} + \ln{a}\ln^2{z}}{z(1 - z)} \, \rmd{d}z  - \int_{\frac{a}{1+a}}^1 \frac{\ln^3{z} - \ln{a}\ln^2{z}}{z(1 - z)} \, \rmd{d}z.
\end{split}
\end{equation}
By applying \eqref{eulan3} to \eqref{geezh} and subsequently utilizing \eqref{newinh1}--\eqref{frsin2} in \eqref{geezh}, we derive
\begin{equation}\label{geezh2}
\begin{split}
&\li_4\left(-\frac{1}{a}\right) + \li_4(-a) + \left(\ln{a}\ln(1 + a) - \ln^2(1+a)\right)\left(\li_2\left(\frac{1}{1+a}\right) + \li_2\left(\frac{a}{1+a}\right)\right)
\\&=  -\frac{7\pi^4}{360} - \frac{\pi^2}{12}\ln^2{a} - \frac{\ln^4{a}}{24} + \frac{\pi^2}{6}\ln{a}\ln(1 + a) - \frac{\pi^2}{6}\ln^2(1 + a) - 2\ln{a}\ln^3(1+ a) 
\\&\quad+ \ln^2{a}\ln^2(1 + a) + \ln^4(1 + a).
\end{split}
\end{equation}
Now, applying Euler's reflection formula \cite[(25.12.6)]{bib23} to the dilogarithms in \eqref{geezh2}, we conclude the proof of Theorem \ref{newthgh}. Additionally, it is worth noting that the inversion formula can also be derived by substituting 4 for $m$ and $-\frac{1}{z}$ for $z$ in Jonqui\`ere's inversion formula \cite[\S1.11.1, pp.~31, (16)]{bib2}. However, this newly presented proof offers the advantage of avoiding logarithms of negative real numbers.
\end{proof}

\subsection{Infinite series involving the Hurwitz zeta function}\label{sec3.2}

In this subsection, we introduce a theorem for double infinite series with symmetric summands and apply it to derive new identities. The motivation for determining closed forms for these series stems from the fact that the harmonic series $\sum_{k=1}^\infty \frac{H_k}{k^3}$, as presented in Lemma \ref{lemma4}, can also be proven as follows. By utilizing the series expression for $\psi_1(k)$, we initially obtain
\[\sum_{k=1}^\infty \frac{\psi_1(k)}{k^2} = \sum_{k=1}^\infty \sum_{j=k}^\infty \frac{1}{k^2 j^2}.\]
Interchanging the order of summation, and then the roles of the dummy variables $j$ and $k$, we obtain
\[\sum_{j=1}^\infty \sum_{k=1}^j \frac{1}{j^2k^2} = \sum_{k=1}^\infty \sum_{j=k}^\infty \frac{1}{j^2k^2} = \sum_{k=1}^\infty \sum_{j=1}^k \frac{1}{j^2k^2}.\]
Consequently, this leads to
\[
2\sum_{k=1}^\infty \sum_{j=k}^\infty \frac{1}{j^2k^2} = \zeta(4) + \sum_{k=1}^\infty \sum_{j=1}^\infty \frac{1}{j^2 k^2} = \zeta(4) + \zeta^2(2).\]
Hence,
\[
\sum_{k=1}^\infty \frac{\psi_1(k)}{k^2} = \frac{\zeta(4) + \zeta^2(2)}{2} = \frac{7\pi^4}{360}.
\]
Utilizing \eqref{harmo1}, we derive
\begin{align*}
\sum_{k=1}^\infty \frac{H_k}{k^3} &=  \zeta(4) + \int_0^1 \frac{\left(\frac{\pi^2}{6} - \ln{t}\ln(1- t) - 2\li_2(t)\right)\ln(1 - t)}{t} \, \rmd{d}t + \frac{1}{2} \int_0^1 \frac{\ln{t} \ln^2(1 -t)}{t}\, \rmd{d}t
\\&= \frac{\pi^4}{90} -  \frac{1}{2} \int_0^1 \frac{\ln{t} \ln^2(1 -t)}{t}\, \rmd{d}t = \frac{\pi^4}{90} + \frac{1}{2} \sum_{k=1}^\infty \frac{1}{k} \int_0^1 t^{k-1} \ln{t} \ln(1 - t) \, \rmd{d}t 
\\&= \frac{\pi^4}{60} + \frac{1}{2}\sum_{k=1}^\infty \frac{H_k}{k^3} - \frac{1}{2}\sum_{k=1}^\infty \frac{\psi_1(k)}{k^2} = \frac{\pi^4}{72}.
\end{align*}

\begin{theorem} Let $f$ be an arbitrary symmetric function of two variables such that  \label{theoremh9}
\begin{enumerate}
\item if $f(k, j) = (-1)^{k+j} g(k, j)$, where $g(k, j)$ is positive, and $\sum_{k=1}^\infty f(k, k)$ converges, or
\item if $f(k, j)$ is positive and $\sum_{j=1}^\infty \sum_{k=1}^\infty f(k, j)$ converges.
\end{enumerate}
Then $\sum_{j=1}^\infty \sum_{k=0}^\infty f(k+j, j)$ is convergent, and 
\begin{equation} \label{bigmastht1}
\sum_{j=1}^\infty \sum_{k=0}^\infty f(k+j, j) = \frac{1}{2}\left(\sum_{j=1}^\infty \sum_{k=1}^\infty f(k, j)  + \sum_{k=1}^\infty f(k, k)\right).
\end{equation}
\end{theorem}

\begin{proof}
Interchanging the order of summation, we have
\begin{equation}\label{symmh0}
\sum_{j=1}^\infty \sum_{k=1}^j  f(k, j) = \sum_{k=1}^\infty \sum_{j=k}^\infty f(k, j).
\end{equation}
Next, we interchange the roles of the dummy variables, resulting in
\begin{equation}\label{symmh1}
\sum_{j=1}^\infty \sum_{k=1}^j  f(k, j) =  \sum_{k=1}^\infty \sum_{j=1}^k  f(j,k).
\end{equation}
Since by hypothesis $f(k, j)$ is symmetric, it follows that $f(j, k) = f(k, j)$, and thus, \eqref{symmh1} can be expressed as
\begin{equation}\label{symmh2}
\sum_{j=1}^\infty \sum_{k=1}^j  f(k, j) =  \sum_{k=1}^\infty \sum_{j=1}^k  f(k, j).
\end{equation}
Combining the two resulting expressions from \eqref{symmh0} and \eqref{symmh2}, we deduce
\begin{equation}\label{symmh3}
\sum_{j=1}^\infty \sum_{k=1}^j  f(k, j) = \frac{1}{2}\left(\sum_{j=1}^\infty \sum_{k=1}^\infty f(k, j)  +  \sum_{k=1}^\infty f(k, k)\right),
\end{equation}
By reindexing the series on the left-hand side of \eqref{symmh3}, we have
\begin{equation}\label{symmh4}
\sum_{j=1}^\infty \sum_{k=1}^\infty  f(k, j) - \sum_{j=1}^\infty \sum_{k=j+1}^\infty  f(k, j) =  \frac{1}{2}\left(\sum_{j=1}^\infty \sum_{k=1}^\infty f(k, j)  +  \sum_{k=1}^\infty f(k, k)\right).
\end{equation}
Shifting the index in the second series on the left-hand side of \eqref{symmh3}, we obtain
\begin{equation}\label{symmh5}
\sum_{j=1}^\infty \sum_{k=1}^\infty  f(k, j) - \sum_{j=1}^\infty \sum_{k=1}^\infty  f(k+j, j) =  \frac{1}{2}\left(\sum_{j=1}^\infty \sum_{k=1}^\infty f(k, j)  + \sum_{k=1}^\infty f(k, k)\right).
\end{equation}
Upon reindexing the second series on the left-hand side of \eqref{symmh5} as
\[\sum_{j=1}^\infty \sum_{k=1}^\infty  f(k+j, j) = \sum_{j=1}^\infty \sum_{k=0}^\infty  f(k+j, j) - \sum_{k=1}^\infty f(k, k),\]
and rearranging, \eqref{bigmastht1} follows. We must now demonstrate the conditions for convergence. To begin with the first condition, let $f(k, j) = (-1)^{k+j} g(k, j)$. It is evident that the symmetricity of $f$ implies the symmetricity of $g$. Consider the following series
\[ \sum_{j=1}^\infty \sum_{k=1}^\infty (-1)^{k+j} g(k, j) = \sum_{j=1}^\infty g(k, k) - 2\sum_{1 \leq j < k < \infty} g(k, j) \]
Since, according to our hypothesis, $g(k, j)$ is positive, we can write
\[\sum_{j=1}^\infty \sum_{k=1}^\infty  f(k, j) = \sum_{j=1}^\infty g(k, k) - 2\sum_{1 \leq j < k < \infty} g(k, j) \leq \sum_{j=1}^\infty g(k, k) = \sum_{j=1}^\infty f(k, k).\]
By applying the comparison test, we can conclude that $\sum_{j=1}^\infty \sum_{k=1}^\infty f(k, j)$ converges if $\sum_{j=1}^\infty f(k, k)$ converges. As the sum of two convergent series is itself convergent, we can now deduce from \eqref{bigmastht1} that $\sum_{j=1}^\infty \sum_{k=0}^\infty f(k+j, j)$ converges if $\sum_{k=1}^\infty  f(k, k)$ converges. Moving on to the second condition, assume $\sum_{j=1}^\infty \sum_{k=1}^\infty f(k, j)$ converges. It is clear that since $f(k, j)$ is positive, we can write
\[\sum_{k=1}^\infty  f(k, k) \leq  \sum_{j=1}^\infty f(k, k) + 2\sum_{1 \leq j < k < \infty} f(k, j) = \sum_{j=1}^\infty \sum_{k=1}^\infty f(k, j).\]
By employing a similar argument, we can conclude that $\sum_{j=1}^\infty \sum_{k=0}^\infty f(k+j, j)$ is indeed convergent, provided $\sum_{j=1}^\infty \sum_{k=1}^\infty f(k, j)$ converges.
\end{proof}

\begin{corollary}Let $f$ be an arbitrary function such that
\begin{enumerate}
\item if $f(k) = (-1)^{k} g(k)$, where $g(k)$ is positive, and $\sum_{k=1}^\infty f^2(k)$ converges, or
\item if $f(k)$ is positive and $\sum_{k=1}^\infty f(k)$ converges.
\end{enumerate}
Then $\sum_{j=1}^\infty \sum_{k=0}^\infty f(k+j)f(j)$ is convergent and \label{coroh12}
\begin{equation}
\sum_{j=1}^\infty \sum_{k=0}^\infty f(k+j)f(j) = \frac{1}{2}\left(\left(\sum_{k=1}^\infty f(k)\right)^2+  \sum_{k=1}^\infty f^2(k)\right).
\end{equation}
\end{corollary}

\begin{proof}  
Taking $f(k, j) = f(k)g(j)$ in Theorem \ref{theoremh9}, the symmetricity of $f(k, j)$ implies $f(k) = g(k)$. As such, the proof of Corollary \ref{coroh12} is complete.
\end{proof}

\begin{remark}
Corollary \ref{coroh12} possesses the remarkable property of transforming a double infinite series into an expression that consists of the sum of the square of an infinite series and another infinite series. It is worth noting that the results established in Theorem \ref{theoremh9} and Corollary \ref{coroh12} have not been previously presented elsewhere in the existing literature.
\end{remark}

\noindent
We apply Corollary \ref{coroh12} in the following theorem.
\begin{theorem}Let $\Re(m) > 1$, $r, s \in \mathbb{C}$, where $r \neq 0$, $rk \neq s$, for any positive integer $k$. Then \label{bigthmh1}
\begin{equation}
\sum_{k=1}^\infty \frac{ \zeta\left(m, \frac{rk-s}{r}\right) }{(rk-s)^m}= \frac{1}{2r^m}\left(\zeta^2\left(m, \frac{r-s}{r}\right) + \zeta\left(2m, \frac{r-s}{r}\right)\right).
\end{equation}
\end{theorem}

\begin{proof}
By setting $f(k) = \frac{1}{(rk-s)^m}$ in Corollary \ref{coroh12}, we conclude the proof of Theorem \ref{bigthmh1}.
\end{proof}

\begin{remark}
Theorem \ref{bigthmh1} does not appear in the DLMF \cite[\S25.11(xi)]{bib23} and Prudnikov's book \cite[pp.~396--397]{Prud}, where series involving the Hurwitz zeta function are discussed.
\end{remark}

\begin{corollary} \label{thisc}
Let $m$ be any positive integer greater than $1$, $r, s \in \mathbb{C}$, where $r \neq 0$, $rk \neq s$, for any positive integer $k$. Then 
\begin{equation}
\sum_{j=1}^{\infty} \frac{\psi_{m-1}\left(\frac{r j-s}{r}\right)}{(r j-s)^m} = \frac{(-1)^m}{2 r^m} \left(\frac{\psi_{m-1}^2\left(\frac{r-s}{r}\right)}{(m-1)!} +\frac{(m-1)!}{(2m-1)!}\psi_{2m-1}\left(\frac{r-s}{r}\right)\right).
\end{equation}
\end{corollary}

\begin{proof}
Theorem \ref{bigthmh1} reduces to Corollary \ref{thisc}, if we consider $m$ as any positive integer greater than $1$.
\end{proof}

\begin{ex}For any positive integer $m > 1$, we have
\begin{equation}\label{tirtheh2eq}
\sum_{j=1}^{\infty} \frac{\psi_{m-1}\left(\frac{4j-1}{4}\right)}{(4j-1)^m} = (-1)^m 2^{-2m-1} \left(\frac{\psi_{m-1}^2\left(\frac{3}{4}\right)}{(m-1)!} +\frac{(m-1)!}{(2m-1)!}\psi_{2m-1}\left(\frac{3}{4}\right)\right).
\end{equation}
For $m= 2, 3, 4, 5$, we have
\begin{align}
&\sum_{j=1}^{\infty} \frac{\psi_{1}\left(\frac{4j-1}{4}\right)}{(4j-1)^2}  = 2\G^2 - \frac{\G \pi^2}{2} + \frac{\pi^4}{32} + \frac{\psi_3\left(\frac{3}{4}\right)}{192}, \\
&\sum_{j=1}^{\infty} \frac{\psi_{2}\left(\frac{4j-1}{4}\right)}{(4j-1)^3}  = -\frac{\pi^6}{64} + \frac{7\pi^3}{8}\zeta(3) - \frac{49\zeta^2(3)}{4} - \frac{\psi_5\left(\frac{3}{4}\right)}{7680}, \\
&\sum_{j=1}^{\infty} \frac{\psi_{3}\left(\frac{4j-1}{4}\right)}{(4j-1)^4}  = \frac{\psi_3^2\left(\frac{3}{4}\right)}{3072} + \frac{\psi_7\left(\frac{3}{4}\right)}{430080},\\
&\sum_{j=1}^{\infty} \frac{\psi_{4}\left(\frac{4j-1}{4}\right)}{(4j-1)^5}  = -\frac{25\pi^{10}}{768} + \frac{155 \pi^5}{8} \zeta(5) - 2883\zeta^2(5) - \frac{\psi_9\left(\frac{3}{4}\right)}{30965760}.
\end{align}
\end{ex}

\begin{corollary}\label{corsha}
Let $m$ be any positive integer greater than $1$. Then 
\begin{equation}
\begin{split}
\sum_{j=1}^{\infty} \frac{\psi_{2m-2}\left(\frac{4j-1}{4}\right)}{(4j-1)^{2m-1}} &=   - \frac{|E_{2m-2}|}{8}\left(1-2^{2m-1}\right) \pi^{2m-1} \zeta(2m-1)  
\\&\quad- \frac{\left(1-2^{2m-1}\right)^2(2m-2)!}{8}\zeta^2(2m-1)- \frac{E^2_{2m-2}}{32(2m-2)!}\pi^{4m-2}  
\\&\quad-\frac{(2m-2)!}{2^{4m-1}(4m-3)!}\psi_{4m-3}\left(\frac{3}{4}\right),
\end{split}
\end{equation}
where $E_m$ are the Euler numbers.
\end{corollary}

\begin{proof}
Olaikhan \cite[\S1.20.6, \S1.20.7, pp.~62--63]{bibash} expressed $\psi_{2a}\left(\frac{3}{4}\right)$ in terms of $\psi_{2a}\left(\frac{1}{4}\right)$, and $\psi_{2a}\left(\frac{1}{4}\right)$ in terms of the Euler numbers $E_a$, with $a$ as a positive integer. By employing both of these expressions, we derive
\begin{equation}\label{psiabd1}
\psi_{2m-2}\left(\frac{3}{4}\right) = 2^{2m-2}\left(\left(1 - 2^{2m-1}\right) (2m-2)! \zeta(2m-1) + \frac{\pi^{2m-1}}{2}|E_{2m-2}|\right).
\end{equation}
By replacing $m$ with $2m-1$ and subsequently substituting \eqref{psiabd1} into \eqref{tirtheh2eq}, we conclude the proof of Corollary \ref{corsha}.
\end{proof}

\begin{theorem}Let $m$ be any positive integer greater than 1, $r, s \in \mathbb{C}$, where $r \neq 0$, $rk \neq s$, for any positive integer $k$. Then \label{bghaft}
\begin{align*}
&\sum_{k=1}^\infty \sum_{p=0}^{m-1} \frac{1}{r^p (rk-s)^{m-p}} \binom{m}{p} \sin\left(\frac{\pi}{2}(m-p)\right) \int_0^\infty \frac{x^{m-p}}{(r^2 x^2 + (rk - s)^2)^m (e^{2\pi x} - 1)} \, \rmd{d}x 
\\&\qquad= \frac{1}{4r^{3m}} \zeta^2\left(m, \frac{r-s}{r}\right) - \frac{1}{2r^{3m} (m -1)} \zeta\left(2m-1, \frac{r-s}{r}\right).
\end{align*}
\end{theorem}

\begin{proof}
Employing Hermite's integral representation \eqref{herm} for $\zeta(s, z)$ in Theorem \ref{bigthmh1}, we have
\begin{equation}\label{sphaf}
\begin{split}
&\sum_{k=1}^\infty \frac{1}{(rk-s)^m}\left(\frac{\left(\frac{rk-s}{r}\right)^{-m}}{2} + \frac{\left(\frac{rk-s}{r}\right)^{1-m}}{m-1}\right)
\\&\qquad=  \frac{1}{2r^m}\zeta\left(2m, \frac{r-s}{r}\right) + \frac{1}{r^m(m-1)}\zeta\left(2m-1, \frac{r-s}{r}\right).
\end{split}
\end{equation}
Define the last integral in \eqref{herm} as
\[f(m, z) = \int_0^\infty \frac{\sin\left(m\arctan\left(x/(z)\right)\right)}{\left(x^2 + z^2\right)^{\frac{m}{2}}\left(\mathrm{e}^{2\pi x} - 1\right)}\, \mathrm{d}x.\]
By De Moivre's theorem \cite[(4.21.34)]{bib23}, we have
\[f(m, z) = \Im \int_0^\infty \frac{\left(\cos\left(\arctan\left(\frac{x}{z}\right)\right) + \mathrm{i}\sin\left(\arctan\left(\frac{x}{z}\right)\right)\right)^m}{\left(x^2 + z^2\right)^{\frac{m}{2}}\left(\mathrm{e}^{2\pi x} - 1\right)}\, \mathrm{d}x,\]
where $\mathrm{i}=\sqrt{-1}$. By the binomial theorem, we have
\[f(m, z) = \Im \int_0^\infty\sum_{p=0}^m \binom{m}{p} \frac{\left(\cos\left(\arctan\left(\frac{x}{z}\right)\right)\right)^p \mathrm{i}^{m-p} \left(\sin\left(\arctan\left(\frac{x}{z}\right)\right)\right)^{m-p}}{\left(x^2 + z^2\right)^{\frac{m}{2}}\left(\mathrm{e}^{2\pi x} - 1\right)}\, \mathrm{d}x.\]
Interchanging summation and integration, we have
\[f(m, z) = \Im \sum_{p=0}^m \binom{m}{p} \int_0^\infty\frac{\left(\cos\left(\arctan\left(\frac{x}{z}\right)\right)\right)^p \mathrm{i}^{m-p} \left(\sin\left(\arctan\left(\frac{x}{z}\right)\right)\right)^{m-p}}{\left(x^2 + z^2\right)^{\frac{m}{2}}\left(\mathrm{e}^{2\pi x} - 1\right)}\, \mathrm{d}x.\]
Using the trigonometric identity $1 + \tan^2\theta = \sec^2\theta$, we deduce
\[\cos\left(\arctan\left(\frac{x}{z}\right)\right)  = \frac{z}{\sqrt{x^2 + z^2}}, \quad \sin\left(\arctan\left(\frac{x}{z}\right)\right) = \frac{x}{\sqrt{x^2 + z^2}}.\]
Therefore, we can express $f(m, z)$ as
\[f(m, z) = \Im \sum_{p=0}^m \binom{m}{p} z^p \mathrm{i}^{m-p}\int_0^\infty \frac{x^{m-p}}{\left(x^2 + z^2\right)^{m} \left(\mathrm{e}^{2\pi x} - 1\right)}\, \mathrm{d}x.\]
Since
\[\Im  \mathrm{i}^{m-p} = \Im e^{\mathrm{i}\frac{\pi}{2}(m-p)} = \Im\left(\cos\left(\frac{\pi}{2}(m-p)\right) +\mathrm{i}\sin\left(\frac{\pi}{2}(m-p)\right)\right) = \sin\left(\frac{\pi}{2}(m-p)\right),\]
we conclude that
\begin{equation}\label{sether1}
f(m, z) = \sum_{p=0}^{m-1} \binom{m}{p} z^p  \sin\left(\frac{\pi}{2}(m-p)\right) \int_0^\infty \frac{x^{m-p}}{\left(x^2 + z^2\right)^{m} \left(\mathrm{e}^{2\pi x} - 1\right)}\, \mathrm{d}x.
\end{equation}
Replacing $z$ in \eqref{sether1} with $\frac{rk-s}{r}$, utilizing the latter and \eqref{sphaf} to derive an integral representation for $\zeta\left(m, \frac{rk-s}{r}\right)$, and further substituting the result in Theorem \ref{bigthmh1}, we conclude the proof of Theorem \ref{bghaft}.
\end{proof}
Since $\sin\left(\frac{\pi}{2}(m - p)\right) = 0, \pm 1$, for integer values of $m$ and $p$, Theorem \ref{bghaft} can be rewritten for odd and even values of $m$, as follows.
\begin{corollary} Let $m$ be any positive integer, $r, s \in \mathbb{C}$, where $r \neq 0$, $rk \neq s$, for any positive integer $k$. Then \label{bghaft1o}
\begin{align*}
&\sum_{k=1}^\infty \sum_{p=0}^{m-1} \frac{(-1)^{m+p-1}}{r^{2p+1} (rk-s)^{2m-2p-1}} \binom{2m}{2p+1}  \int_0^\infty \frac{x^{2m-2p-1}}{(r^2 x^2 + (rk - s)^2)^{2m} (e^{2\pi x} - 1)} \, \rmd{d}x 
\\&\qquad= \frac{1}{4r^{6m}} \zeta^2\left(2m, \frac{r-s}{r}\right) - \frac{1}{2r^{6m} (2m -1)} \zeta\left(4m-1, \frac{r-s}{r}\right).
\end{align*}
\end{corollary}
\begin{corollary} Let $m$ be any positive integer, $r, s \in \mathbb{C}$, where $r \neq 0$, $rk \neq s$, for any positive integer $k$. Then \label{bghaft2o}
\begin{align*}
&\sum_{k=1}^\infty \sum_{p=0}^{m} \frac{(-1)^{m+p}}{r^{2p} (rk-s)^{2m-2p+1}} \binom{2m+1}{2p} \int_0^\infty \frac{x^{2m-2p+1}}{(r^2 x^2 + (rk - s)^2)^{2m+1} (e^{2\pi x} - 1)} \, \rmd{d}x 
\\&\qquad= \frac{1}{4r^{6m+3}} \zeta^2\left(2m+1, \frac{r-s}{r}\right) - \frac{1}{4 r^{6m+3} m} \zeta\left(4m+1, \frac{r-s}{r}\right).
\end{align*}
\end{corollary}

\noindent
Theorem \ref{bghaft} reduces directly to an expression in the Riemann zeta function if we set $r=1$, $s=0$ and $r=2$, $s=1$. This yields the following corollaries.

\begin{corollary}Let $m$ be any positive integer greater than 1. Then \label{corol1h1h}
\begin{align*}
&\sum_{k=1}^\infty \sum_{p=0}^{m-1} \frac{1}{k^{m-p}} \binom{m}{p} \sin\left(\frac{\pi}{2}(m-p)\right)  \int_0^\infty \frac{x^{m-p}}{(x^2 + k^2)^m (e^{2\pi x} - 1)} \, \rmd{d}x 
\\&\qquad= \frac{1}{4} \zeta^2\left(m\right) - \frac{1}{2(m -1)} \zeta\left(2m-1\right).
\end{align*}
\end{corollary}

\begin{corollary}Let $m$ be any positive integer greater than 1. Then \label{corol1h}
\begin{align*}
&\sum_{k=1}^\infty \sum_{p=0}^{m-1} \frac{1}{2^p (2k-1)^{m-p}} \binom{m}{p} \sin\left(\frac{\pi}{2}(m-p)\right) \int_0^\infty \frac{x^{m-p}}{(4 x^2 + (2k - 1)^2)^m (e^{2\pi x} - 1)} \, \rmd{d}x 
\\&\qquad= 2^{-3 m-2} \left(2^m-1\right)^2 \zeta^2(m) - \frac{2^{-3 m-1} \left(2^{2 m-1}-1\right)}{m-1}  \zeta (2 m-1).
\end{align*}
\end{corollary}

\noindent
As in Corollaries \ref{bghaft1o} and \ref{bghaft2o}, we provide expressions for Corollaries \ref{corol1h1h} and \ref{corol1h} for odd and even values of $m$.

\begin{corollary}Let $m$ be any positive integer. Then \label{corol1}
\begin{align*}
&\sum_{k=1}^\infty \sum_{p=0}^{m-1} \frac{(-1)^{m+p-1}}{k^{2m-2p-1}} \binom{2m}{2p+1}  \int_0^\infty \frac{x^{2m-2p-1}}{(x^2 + k^2)^{2m} (e^{2\pi x} - 1)} \, \rmd{d}x 
\\&\qquad= \frac{2^{4m-4}}{(2m)!^2} B_{2m}^2 \pi^{4m}  - \frac{1}{2 (2m -1)} \zeta(4m-1),
\end{align*}
where $B_m$ are the Bernoulli numbers.
\end{corollary}

\begin{ex}For $m=1, 2, 3$, we have
\[\sum_{k=1}^\infty  \frac{1}{k}  \int_0^\infty \frac{x}{(x^2 + k^2)^2 (e^{2\pi x} - 1)} \, \rmd{d}x =  \frac{\pi^4}{288} - \frac{\zeta(3)}{4},\]

\[\sum_{k=1}^\infty  \frac{1}{k} \int_0^\infty\frac{x}{\left(x^2 + k^2\right)^4\left(e^{2 \pi  x}-1\right)} \, \rmd{d}x - \sum_{k=1}^\infty  \frac{1}{k^3} \int_0^\infty \frac{x^{3}}{\left(x^2 + k^2\right)^4\left(e^{2 \pi  x}-1\right)} \, \rmd{d}x = \frac{\pi^8}{129600}-\frac{\zeta(7)}{24},\]

\begin{align*}
& 6 \sum_{k=1}^\infty  \frac{1}{k} \int_0^\infty\frac{x}{\left(x^2 + k^2\right)^6\left(e^{2 \pi  x}-1\right)} \, \rmd{d}x - 20\sum_{k=1}^\infty  \frac{1}{k^3} \int_0^\infty \frac{x^{3}}{\left(x^2 + k^2\right)^6\left(e^{2 \pi  x}-1\right)} \, \rmd{d}x 
\\&\qquad + 6\sum_{k=1}^\infty  \frac{1}{k^5} \int_0^\infty \frac{x^5}{\left(x^2 + k^2\right)^6 \left(e^{2 \pi  x}-1\right)} \, \rmd{d}x = \frac{\pi ^{12}}{3572100}-\frac{\zeta(11)}{10}.
\end{align*}
\end{ex}

\begin{corollary}Let $m$ be any positive integer. Then \label{corol2}
\begin{align*}
&\sum_{k=1}^\infty \sum_{p=0}^{m} \frac{(-1)^{m+p}}{k^{2m-2p+1}} \binom{2m+1}{2p} \int_0^\infty \frac{x^{2m-2p+1}}{(x^2 + k^2)^{2m+1} (e^{2\pi x} - 1)} \, \rmd{d}x 
\\&\qquad= \frac{1}{4} \zeta^2(2m+1) - \frac{1}{4 m} \zeta(4m+1).
\end{align*}
\end{corollary}

\begin{ex}For $m=1, 2, 3$, we have
\[3 \sum_{k=1}^\infty  \frac{1}{k} \int_0^\infty\frac{x}{\left(x^2 + k^2\right)^3\left(e^{2 \pi  x}-1\right)} \, \rmd{d}x -\sum_{k=1}^\infty  \frac{1}{k^3} \int_0^\infty \frac{x^{3}}{\left(x^2 + k^2\right)^3\left(e^{2 \pi  x}-1\right)} \, \rmd{d}x = \frac{\zeta^2(3)}{4}-\frac{\zeta(5)}{4},\]

\begin{align*}
& 5\sum_{k=1}^\infty  \frac{1}{k} \int_0^\infty\frac{x}{\left(x^2 + k^2\right)^5\left(e^{2 \pi  x}-1\right)} \, \rmd{d}x  - 10\sum_{k=1}^\infty  \frac{1}{k^3} \int_0^\infty \frac{x^{3}}{\left(x^2 + k^2\right)^5\left(e^{2 \pi  x}-1\right)} \, \rmd{d}x 
\\&\qquad + \sum_{k=1}^\infty  \frac{1}{k^5} \int_0^\infty \frac{x^5}{\left(x^2 + k^2\right)^5\left(e^{2 \pi  x}-1\right)} \, \rmd{d}x = \frac{\zeta^2(5)}{4}-\frac{\zeta(9)}{8},
\end{align*}

\begin{align*}
&  7 \sum_{k=1}^\infty  \frac{1}{k} \int_0^\infty\frac{x}{\left(x^2 + k^2\right)^7\left(e^{2 \pi  x}-1\right)} \, \rmd{d}x  - 35\sum_{k=1}^\infty  \frac{1}{k^3} \int_0^\infty \frac{x^{3}}{\left(x^2 + k^2\right)^7\left(e^{2 \pi  x}-1\right)} \, \rmd{d}x 
\\&\qquad + 21\sum_{k=1}^\infty  \frac{1}{k^5} \int_0^\infty \frac{x^5}{\left(x^2 + k^2\right)^7\left(e^{2 \pi  x}-1\right)} \, \rmd{d}x - \sum_{k=1}^\infty  \frac{1}{k^7} \int_0^\infty \frac{x^7}{\left(x^2 + k^2\right)^7\left(e^{2 \pi  x}-1\right)} \, \rmd{d}x
\\&\qquad= \frac{\zeta^2(7)}{4}-\frac{\zeta(13)}{12}.
\end{align*}
\end{ex}

\begin{corollary}Let $m$ be any positive integer. Then \label{corol3}
\begin{align*}
&\sum_{k=1}^\infty \sum_{p=0}^{m-1} \frac{(-1)^{m+p-1}}{2^{2p+1} (2k-1)^{2m-2p-1}} \binom{2m}{2p+1} \int_0^\infty \frac{x^{2m-2p-1}}{(4 x^2 + (2k - 1)^2)^{2m} (e^{2\pi x} - 1)} \, \rmd{d}x 
\\&\qquad= \frac{2^{-2m-4} \left(2^{2m}-1\right)^2}{(2m)!^2}  B^2_{2m} \pi^{4m} - \frac{2^{-6m-1} \left(2^{4m-1}-1\right)}{2m-1}  \zeta (4m-1),
\end{align*}
where $B_m$ are the Bernoulli numbers.
\end{corollary}

\begin{ex}For $m=1, 2, 3$, we have
\[\sum_{k=1}^\infty  \frac{1}{2k-1}  \int_0^\infty \frac{x}{(4 x^2 + (2k - 1)^2)^2 (e^{2\pi x} - 1)} \, \rmd{d}x = \frac{\pi ^4}{1024}-\frac{7 \zeta (3)}{128},\]

\begin{align*}
& \frac{1}{2} \sum_{k=1}^\infty  \frac{1}{2 k-1} \int_0^\infty\frac{x}{\left(4x^2 + (2 k-1)^2\right)^4\left(e^{2 \pi  x}-1\right)} \, \rmd{d}x 
\\&\qquad - 2\sum_{k=1}^\infty  \frac{1}{(2 k-1)^3} \int_0^\infty \frac{x^{3}}{\left(4x^2 + (2 k-1)^2\right)^4\left(e^{2 \pi  x}-1\right)} \, \rmd{d}x = \frac{\pi ^8}{589824} - \frac{127 \zeta (7)}{24576},
\end{align*}

\begin{align*}
& \frac{3}{16} \sum_{k=1}^\infty  \frac{1}{2 k-1} \int_0^\infty\frac{x}{\left(4x^2 + (2 k-1)^2\right)^6\left(e^{2 \pi  x}-1\right)} \, \rmd{d}x 
\\&\qquad - \frac{5}{2}\sum_{k=1}^\infty  \frac{1}{(2 k-1)^3} \int_0^\infty \frac{x^{3}}{\left(4x^2 + (2 k-1)^2\right)^6\left(e^{2 \pi  x}-1\right)} \, \rmd{d}x 
\\&\qquad + 3\sum_{k=1}^\infty  \frac{1}{(2 k-1)^5} \int_0^\infty \frac{x^5}{\left(4x^2 + (2 k-1)^2\right)^6 \left(e^{2 \pi  x}-1\right)} \, \rmd{d}x = \frac{\pi ^{12}}{235929600}-\frac{2047 \zeta (11)}{2621440}.
\end{align*}
\end{ex}

\begin{corollary}Let $m$ be any positive integer. Then \label{corol4}
\begin{align*}
&\sum_{k=1}^\infty \sum_{p=0}^{m} \frac{ (-1)^{m+p}}{2^{2p} (2k-1)^{2m-2p+1}} \binom{2m+1}{2p} \int_0^\infty \frac{x^{2m-2p+1}}{(4 x^2 + (2k - 1)^2)^{2m+1} (e^{2\pi x} - 1)} \, \rmd{d}x 
\\&\qquad= 2^{-6m-5} \left(2^{2m+1}-1\right)^2 \zeta^2(2m+1) - \frac{2^{-6m-4} \left(2^{4m+1}-1\right)}{2m}  \zeta (4m+1).
\end{align*}
\end{corollary}

\begin{ex} For $m=1, 2, 3$, we have
\begin{align*}
&\frac{3}{4} \sum_{k=1}^\infty  \frac{1}{2 k-1} \int_0^\infty\frac{x}{\left(4x^2 + (2 k-1)^2\right)^3\left(e^{2 \pi  x}-1\right)} \, \rmd{d}x 
\\&\qquad-\sum_{k=1}^\infty  \frac{1}{(2 k-1)^3} \int_0^\infty \frac{x^{3}}{\left(4x^2 + (2 k-1)^2\right)^3\left(e^{2 \pi  x}-1\right)} \, \rmd{d}x = \frac{49 \zeta^2(3)}{2048}-\frac{31 \zeta (5)}{2048},
\end{align*}

\begin{align*}
& \frac{5}{16} \sum_{k=1}^\infty  \frac{1}{2 k-1} \int_0^\infty\frac{x}{\left(4x^2 + (2 k-1)^2\right)^5\left(e^{2 \pi  x}-1\right)} \, \rmd{d}x 
\\&\qquad - \frac{5}{2}\sum_{k=1}^\infty  \frac{1}{(2 k-1)^3} \int_0^\infty \frac{x^{3}}{\left(4x^2 + (2 k-1)^2\right)^5\left(e^{2 \pi  x}-1\right)} \, \rmd{d}x 
\\&\qquad + \sum_{k=1}^\infty  \frac{1}{(2 k-1)^5} \int_0^\infty \frac{x^5}{\left(4x^2 + (2 k-1)^2\right)^5\left(e^{2 \pi  x}-1\right)} \, \rmd{d}x = \frac{961 \zeta^2(5)}{131072}-\frac{511 \zeta(9)}{262144},
\end{align*}

\begin{align*}
& \frac{7}{64} \sum_{k=1}^\infty  \frac{1}{2 k-1} \int_0^\infty\frac{x}{\left(4x^2 + (2 k-1)^2\right)^7\left(e^{2 \pi  x}-1\right)} \, \rmd{d}x 
\\&\qquad - \frac{35}{16}\sum_{k=1}^\infty  \frac{1}{(2 k-1)^3} \int_0^\infty \frac{x^{3}}{\left(4x^2 + (2 k-1)^2\right)^7\left(e^{2 \pi  x}-1\right)} \, \rmd{d}x 
\\&\qquad + \frac{21}{4}\sum_{k=1}^\infty  \frac{1}{(2 k-1)^5} \int_0^\infty \frac{x^5}{\left(4x^2 + (2 k-1)^2\right)^7\left(e^{2 \pi  x}-1\right)} \, \rmd{d}x 
\\&\qquad - \sum_{k=1}^\infty  \frac{1}{(2 k-1)^7} \int_0^\infty \frac{x^7}{\left(4x^2 + (2 k-1)^2\right)^7\left(e^{2 \pi  x}-1\right)} \, \rmd{d}x= \frac{16129 \zeta^2(7)}{8388608}-\frac{8191 \zeta(13)}{25165824}.
\end{align*}
\end{ex}

\section{Conclusion}
In this work, we have introduced several new theorems that provide closed forms for generalized integrals and series. Furthermore, we have demonstrated the application of our double infinite series transformation formula in deriving new identities, which have been expressed using well-known numbers such as the Euler and Bernoulli numbers. One that we find interesting are the simplest cases of Theorem \ref{bghaft}:
\begin{align*}
&\int_0^\infty \frac{x}{(x^2 + 1)^2 (e^{2\pi x} - 1)} \, \rmd{d}x + \frac{1}{2} \int_0^\infty \frac{x}{(x^2 + 4)^2 (e^{2\pi x} - 1)} \, \rmd{d}x 
\\&\qquad+  \frac{1}{3} \int_0^\infty \frac{x}{(x^2 + 9)^2 (e^{2\pi x} - 1)} \, \rmd{d}x  +  \frac{1}{4} \int_0^\infty \frac{x}{(x^2 + 16)^2 (e^{2\pi x} - 1)} \, \rmd{d}x  + \cdots 
\\&=  \frac{\pi^4}{288} - \frac{\zeta(3)}{4},
\end{align*}
\begin{align*}
&\int_0^\infty \frac{x}{(4 x^2 + 1)^2 (e^{2\pi x} - 1)} \, \rmd{d}x + \frac{1}{3} \int_0^\infty \frac{x}{(4 x^2 + 9)^2 (e^{2\pi x} - 1)} \, \rmd{d}x 
\\&\qquad+  \frac{1}{5} \int_0^\infty \frac{x}{(4 x^2 + 25)^2 (e^{2\pi x} - 1)} \, \rmd{d}x  +  \frac{1}{7} \int_0^\infty \frac{x}{(4 x^2 + 49)^2 (e^{2\pi x} - 1)} \, \rmd{d}x  + \cdots 
\\&= \frac{\pi ^4}{1024}-\frac{7 \zeta (3)}{128}.
\end{align*}
Interested readers can further explore the transformation formulas outlined in Theorem \ref{theoremh9} and Corollary \ref{coroh12} to potentially uncover additional new results.

\section*{Acknowledgment}
I would like to express my heartfelt gratitude to the Spirit of Ramanujan (SOR) STEM Talent Initiative, directed by Ken Ono, for providing the essential computational tools that facilitated the verification and validation of my results. My sincere appreciation extends to esteemed researchers Ali S.~Olaikhan and Cornel I. V\u{a}lean for their invaluable insights. I also thank Mr.~G.~S.~Lawal, Dr.~H.~P.~Adeyemo, Dr.~D.~A.~Dikko, and Dr.~I.~Adinya from the Department of Mathematics at the University of Ibadan for their unwavering support and encouragement.
\section*{Funding}
The author did not receive funding from any organization for the submitted work.
\bibliographystyle{unsrt}
\bibliography{Research_Dilog}

\begin{thebibliography}{1}

\bibitem{bib23}
{\it NIST Digital Library of Mathematical Functions}.
\newblock \href{https://dlmf.nist.gov/}{\bf\tt\normalsize
  https://dlmf.nist.gov/}, Release 1.1.10 of 2023-06-15.
\newblock F.~W.~J. Olver, A.~B. {Olde Daalhuis}, D.~W. Lozier, B.~I. Schneider,
  R.~F. Boisvert, C.~W. Clark, B.~R. Miller, B.~V. Saunders, H.~S. Cohl, and
  M.~A. McClain, eds.

\bibitem{bib2}
{H. Bateman and Erdélyi}.
\newblock Higher {T}ranscendental {F}unctions.
\newblock {\em McGraw-{H}ill {B}ook, {N}ew {Y}ork}, 1:317, 1955,
  \url{https://doi.org/10.1038/175317a0}.

\bibitem{bib9}
{C. I. V\u{a}lean}.
\newblock {\em \normalfont {M}ore ({A}lmost) {I}mpossible {I}ntegrals, {S}ums,
  and {S}eries}.
\newblock Problem {B}ooks in {M}athematics. Springer, Cham, 2023,
  \url{https://doi.org/10.1007/978-3-031-21262-8}.

\bibitem{bib32}
{C. I. V\u{a}lean}.
\newblock {\em \normalfont({A}lmost) {I}mpossible {I}ntegrals, {S}ums, and
  {S}eries}.
\newblock Problem {B}ooks in {M}athematics. Springer, Cham, 2019,
  \url{https://doi.org/10.1007/978-3-030-02462-8}.

\bibitem{bibch}
{E.~Hairer, G.~Wanner}.
\newblock Analysis by {I}ts {H}istory.
\newblock {\em Springer, New York}, pages 197--199, chapter III, 1996.

\bibitem{Prud}
{A.~P.~Prudnikov, Yu.~A.~Brychkov and O.~I.~Marichev}.
\newblock {Integrals and Series: More Special Functions}.
\newblock {\em {Gordon and Breach Science Publishers, New York}}, 3:800, 1990.

\bibitem{bibash}
{A. S. Olaikhan}.
\newblock {\em \normalfont {A}n {I}ntroduction to the {H}armonic {S}eries and
  {L}ogarithmic {I}ntegrals}.
\newblock {F}or {H}igh {S}chool {S}tudents up to {R}esearchers. Second edition,
  2023, ISBN: 978-1-7367360-3-6.

\end{thebibliography}
\end{document}